\documentclass[12pt]{article}

\usepackage{amsmath}
\usepackage{amssymb}
\usepackage{amscd}
\usepackage{epsfig}

 \newcommand{\qed}{\hfill $\square$}
 \newcommand{\bee}{\begin{equation}}
 \newcommand{\eee}{\end{equation}}
  
 \newcommand{\Lb}{\mbox {\boldmath ${\Lambda}$}}  
 \newcommand{\Gb}{\mbox {\boldmath ${\Gamma}$}} 
\newcommand{\Gbt}{\mbox{\tiny\boldmath ${\Gamma}$}}
 \newcommand{\Lbs}{\mbox{\scriptsize\boldmath ${\Lambda}$}}
 \newcommand{\Lbt}{\mbox{\tiny\boldmath ${\Lambda}$}}
 \newcommand{\Pb}{\mbox {\bf P}}
 \newcommand{\Pbs}{\mbox {\scriptsize{\bf P}}}
 \newcommand{\Qb}{\mbox {\bf Q}}
 \newcommand{\Qbs}{\mbox {\scriptsize{\bf Q}}}
 
\newcommand{\diam}{\mbox{\rm diam}}
\newcommand{\be}{\begin{eqnarray}}
\newcommand{\ee}{\end{eqnarray}}
\newcommand{\supp}{\mbox{\rm supp}}
\newcommand{\suppm}{\mbox{\rm \scriptsize supp}}

\newcommand{\freq}{\mbox{\rm freq}}
\newcommand{\Vol}{\mbox{\rm Vol}}
\newcommand{\es}{\emptyset}

\newcommand{\e}{{\varepsilon}}

\newcommand{\gam}{{\gamma}}

\newcommand{\R}{{\mathbb R}}

\newcommand{\Z}{{\mathbb Z}}
\newcommand{\C}{{\mathbb C}}

\newcommand{\Hk}{{\mathcal H}}

\newcommand{\Qk}{{\mathcal Q}}

\newcommand{\bd}{\partial}

\newcommand{\Cant}{{\mathcal C}}

\newcommand{\Ck}{{\mathcal C}}

\newcommand{\dist}{\mbox{\rm dist}}

\newcommand{\Lam}{{\Lambda}}

\newcommand{\om}{\omega}
\newcommand{\mmin}{\rm min}
\newcommand{\mmax}{\rm max}

\newcommand{\nub}{{\nu_{_{\Lbs}}}}
\newcommand{\nubp}{{\nu_{_{\Lbs'}}}}
\newcommand{\cheiv}{{\chi_{_{{\scriptsize {\bf E}_i},V}}}}
\newcommand{\chipv}{{\chi_{_{{\scriptsize\bf P},V}}}}

 \newtheorem{theorem}{Theorem}[section]
 \newtheorem{lemma}[theorem]{Lemma}
 
 \newtheorem{cor}[theorem]{Corollary}

 \newtheorem{defi}[theorem]{Definition}

\numberwithin{equation}{section}

\begin{document}
\bibliographystyle{unset}
 
\title{{\sc Pure Point Dynamical and Diffraction Spectra}}
\date \today 
\maketitle

\vspace{7mm}
 \centerline {
{\sc Jeong-Yup  Lee $^{\,\rm a}$},
{\sc Robert V.\ Moody $^{\,\rm a}$},
{\sc Boris Solomyak $^{\,\rm b}$\footnote{BS acknowledges support from
NSF grants DMS 9800786 and DMS 0099814.}
}}

\vspace{7mm}

{\small
\hspace*{4em}
  
\smallskip
\hspace*{4em}
a:  Dept. of Mathematical and Statistical Sciences, University of Alberta, \\
\hspace*{4em}
\hspace*{2.5em} Edmonton, Alberta T6G 2G1, Canada 

\smallskip
\hspace*{4em}
b:  Department of Mathematics, University of Washington\\
\hspace*{4em}
\hspace*{2.5em} Seattle, WA 98195, USA.}

\abstract
We show that for multi-colored Delone point sets with finite local complexity and uniform cluster frequencies the notions of pure point diffraction and pure point dynamical spectrum are equivalent.

\section{Introduction}
The notion of pure pointedness appears in the theory of aperiodic systems
in two different forms: pure point dynamical spectrum and pure point 
diffraction spectrum. The objective of this paper is to show that these two widely used notions are equivalent under a type of statistical condition known as the existence of uniform cluster frequencies. 

The basic objects of study here are Delone point sets of $\R^d$. The points
of these sets are permitted to be multi-colored, the colors coming from
a finite set of colors. We also assume that our point sets $\Lb$ have the 
property of finite local complexity (FLC), which informally means that
there are only finitely many translational classes of clusters of $\Lb$ 
with any given size. Under these circumstances, the orbit of $\Lb$ under translation gives rise, via completion in the standard Radin-Wolff
type topology, to a compact space $X_{\Lbs}$. With the obvious action
of $\R^d$, we obtain a dynamical system $(X_{\Lbs}, \R^d)$. 

The {\em dynamical spectrum} refers to the spectrum of this 
dynamical system, that is to say, the spectrum of the unitary 
operators $U_x$ arising from the translational action on the 
space of $L^2$-functions on $X_{\Lbs}$.

On the other hand, the {\em diffraction spectrum} (which is the idealized 
mathematical interpretation of the diffraction pattern of a physical 
experiment) is obtained by first assigning weights to the various colors of 
the multiset and then determining the autocorrelation, if it exists,
of this weighted multiset. The Fourier transform of the autocorrelation is the {\em diffraction measure} whose pure pointedness is the question.

There is a well known argument of S. Dworkin (\cite{Dwor}, \cite{Hof}) 
that shows how to deduce pure pointedness of the diffractive spectrum 
from pure pointedness of the dynamical system. Our main result
(Theorem \ref{thm-equiv}) shows that, under the additional assumption 
that $\Lb$ has uniform cluster frequencies (or equivalently, that the 
dynamical system $X_{\Lbs}$ is uniquely ergodic), the process 
can be reversed, 
so in fact the two notions of pure pointedness are equivalent. 

The present understanding of diffractive point sets is very limited.
One of the important consequences of this result is that it allows the introduction of powerful spectral theorems in the study of such sets. Our forthcoming paper \cite{LMS} on diffractive substitution systems makes extensive use of this connection. 

In the last section we discuss what can be salvaged when there are no
uniform cluster frequencies. Then the equivalence of pure point
dynamical and pure point diffraction spectra still holds --- perhaps, not
for the original Delone set $\Lb$, but for almost every Delone set 
in $X_{\Lbs}$, with respect to an ergodic invariant measure.

The prototype of the dynamical system $(X_{\Lbs},\R^d)$
is a symbolic dynamical system, that is, the $\Z$-action by shifts
on a space of bi-infinite sequences.
In the symbolic setting, the equivalence of
pure point dynamical and diffraction spectra has been established
by Queffelec \cite[Prop.\ IV.21]{Queff}, 
and our proof is largely a generalization of her argument.

When the dynamical spectrum is not pure point, its relation to the
diffraction spectrum is not completely understood. It follows from
\cite{Dwor} that the latter is essentially a ``part'' of the former.
So, for instance, if the dynamical spectrum is pure 
singular/absolutely continuous, then the diffraction spectrum is
pure singular/absolutely continuous (apart from the trivial constant
eigenfunction which corresponds to a delta function at $0$).
However, the other direction is more delicate: Van Enter and Mi\c{e}kisz
\cite{VEM} have pointed out that, in the case of mixed spectrum,
the non-trivial pure point component may be ``lost'' when passing from
dynamical spectrum to diffraction spectrum. 

The presentation below contains a number of results that are essentially
well-known, though not always quite in the form needed here. 
For the convenience
 of the reader we have attempted to make the paper largely self-contained. 


\section{Multisets, dynamical systems, and uniform cluster frequencies}

\noindent
A {\em multiset} or {\em $m$-multiset} in $\R^d$ is a 
subset $\Lb = \Lam_1 \times \dots \times \Lam_m 
\subset \R^d \times \dots \times \R^d$ \; ($m$ copies)
where $\Lam_i \subset \R^d$. We also write 
$\Lb = (\Lam_1, \dots, \Lam_m) = (\Lam_i)_{i\le m}$.
We say that $\Lb=(\Lambda_i)_{i\le m}$ is a {\em Delone multiset} in $\R^d$ if
each $\Lambda_i$ is Delone and $\supp(\Lb):=\bigcup_{i=1}^m \Lambda_i 
\subset \R^d$ is Delone.

Although $\Lb$ is a product of sets, it is convenient to think
of it as a set with types or colors, $i$ being the
color of points in $\Lambda_i$. A {\em cluster} of $\Lb$ is, 
by definition,
a family $\Pb = (P_i)_{i\le m}$ where $P_i \subset \Lambda_i$ is 
finite for all $i\le m$.
Many of the clusters that we consider have the form
$A\cap \Lb := (A\cap \Lambda_i)_{i\le m}$, for a bounded set
 $A\subset \R^d$.
There is a natural
translation $\R^d$-action on the set of Delone multisets and their clusters
in $\R^d$. 
The translate of a cluster $\mbox{\bf P}$ by $x \in \R^d$ is 
$x + \Pb = (x + P_i)_{i\le m}$.
We say that two clusters $\Pb$ and $\Pb'$ are {\em translationally equivalent}
if $\Pb=x+\Pb'$ for some $x \in \R^d$. 

We write $B_R(y)$ for the {\em closed} ball of radius $R$
centered at $y$ and use also $B_R$ for $B_R(0)$.

\begin{defi} \label{def-flc}
The Delone multiset $\Lb$ has {\em finite local complexity
(FLC)} if for every $R>0$ there exists a finite set $Y\subset \supp(\Lb)=
\bigcup_{i=1}^m \Lam_i$  such that 
$$
\forall x\in \supp(\Lb),\ \exists\, y\in Y:\ 
B_R(x) \cap \Lb = (B_R(y) \cap \Lb) + (x-y).
$$
In plain language, for each radius $R > 0$ there are only finitely many 
translational classes of clusters whose support lies in some ball of 
radius $R$.
\end{defi}

In this paper we will usually assume that our Delone multisets 
have FLC.

\medskip

Let $\Lb$ be a Delone multiset and $X$ be the collection of all Delone multisets each of whose
clusters is a translate of a $\Lb$-cluster. We introduce a metric
on Delone multisets in a simple variation of the standard way :
 for Delone multisets $\Lb_1$, $\Lb_2 \in X$,
\be \label{metric-multisets}
d(\Lb_1,\Lb_2) := \min\{\tilde{d}(\Lb_1,\Lb_2), 2^{-1/2}\}\, ,
\ee
where
\be 
\tilde{d}(\Lb_1,\Lb_2)
&=&\mbox{inf} \{ \e > 0 : \exists~ x,y \in B_{\e}(0), \nonumber \\ \nonumber
&  & ~~~~~~~~~~ B_{1/{\e}}(0) \cap (-x + \Lb_1) = B_{1/{\e}}(0) 
\cap (-y + \Lb_2) \}\,. 
\ee

Let us indicate why this is a metric. Clearly, the only issue is the
triangle inequality. Suppose that $d(\Lb_1,\Lb_2) \le \e_1,\
d(\Lb_2,\Lb_3) \le \e_2$; we want to show that $d(\Lb_1,\Lb_3) \le \e_1+\e_2$.
We can assume that 
$\e_1,\e_2< 2^{-1/2}$, otherwise the claim is obvious. Then
$$
(-x_1+ \Lb_1) \cap B_{1/\e_1}(0) = (-x_2 + \Lb_2) \cap B_{1/\e_1}(0)\ \ 
\mbox{for some}\ x_1,x_2\in B_{\e_1}(0),
$$
$$
(-x'_2 + \Lb_2) \cap B_{1/\e_2}(0) = (-x'_3 + \Lb_3) \cap B_{1/\e_2}(0)\ \ 
\mbox{for some}\ x'_2, x'_3\in B_{\e_2}(0).
$$
It follows that
$$
(-x_1-x'_2 + \Lb_1) \cap B_{1/\e_1}(-x_2') 
= (-x_2-x'_2 +\Lb_2) \cap B_{1/\e_1}(-x_2').
$$
Since $B_{1/\e_1}(-x_2')\supset B_{(1/\e_1)-\e_2}(0)$, this implies
\be \label{metric1}
(-x_1-x'_2 + \Lb_1) \cap B_{(1/\e_1)-\e_2}(0)
= (-x_2-x'_2 +\Lb_2) \cap B_{(1/\e_1) -\e_2}(0).
\ee
Similarly,
\be \label{metric2}
(-x_2-x'_2 +\Lb_2) \cap B_{(1/\e_2)-\e_1}(0)
= (-x_2-x'_3 + \Lb_3) \cap B_{(1/\e_2)-\e_1}(0).
\ee
A simple computation shows that $\frac{1}{\e_1}-\e_2 \ge \frac{1}{\e_1+\e_2}$
and $\frac{1}{\e_2}-\e_1 \ge \frac{1}{\e_1+\e_2}$ when $\e_1,\e_2< 2^{-1/2}$,
so by (\ref{metric1}) and (\ref{metric2}),
$$
(-x_1-x'_2 + \Lb_1) \cap B_{1/(\e_1+\e_2)}(0)
= (-x_2-x'_3 + \Lb_3) \cap B_{1/(\e_1+\e_2)}(0),
$$
hence $d(\Lb_1,\Lb_3) \le \e_1+\e_2$.

\medskip

We define $X_{\Lbs} := \overline{\{-h + \Lb : h \in \R^d \}}$ with the metric 
$d$. 
In spite of the special role 
played by $0$ in the definition of $d$, any other point of $\R^d$ may 
be used as a reference point, leading to an equivalent metric and more
 importantly the same topology on $X_{\Lbs}$.
The following lemma is standard.
\begin{lemma}\label{lem-comp} {\em (\cite{RW}, \cite{martin})}
If $\Lb$ has FLC, then the metric space $X_{\Lbs}$ is compact.
\end{lemma}
 
The group $\R^d$ acts on $X_{\Lbs}$ by translations which are obviously
homeomorphisms, and we get a topological dynamical system $(X_{\Lbs},\R^d)$.

\begin{defi} \label{def-cyl}
Let $\Pb$ be a cluster of $\Lb$ or some translate of $\Lb$, and let 
$V \subset \R^d$ be a Borel set. 
Define the {\em cylinder set} $X_{\Pbs,V}\subset X_{\Lbs}$ by
\[X_{\Pbs,V} := \{ \Lb' \in X_{\Lbs}:\, -g+\Pb \subset \Lb' \ \mbox{for some}\ 
g\in V\}.\]
\end{defi}

Let $\eta(\Lb) > 0$ be chosen so that every ball of radius 
$\frac{\eta(\Lb)}{2}$ 
contains at most one point of $\supp(\Lb)$, and let $b(\Lb)>0$ be such that 
every ball of radius $\frac{b(\Lb)}{2}$ 
contains at least a point in $\supp(\Lb)$.
These exist by the Delone set property.

The following technical result will be quite useful.

\begin{lemma} \label{lem-tech} Let $\Lb$ be a Delone multiset with FLC.
For any $R\ge \frac{b(\Lb)}{2}$ and $0 < \delta < \eta(\Lb)$, there exist 
Delone multisets $\Gb_j \in X_{\Lbs}$ and Borel sets $V_j$ with  
$\diam(V_j) <\delta$, 
$\Vol(\bd V_j) = 0$, $1 \le j \le N$, such that 
$$
X_{\Lbs} = \bigcup_{j=1}^N X_{\Pbs_j,V_j}
$$
is a disjoint union, where $\Pb_j = B_R(0) \cap \Gb_j$. 
\end{lemma}

\noindent {\em Proof.} For any $R \ge \frac{b(\Lb)}{2}$ consider the clusters 
$\{B_R(0) \cap \Gb :\ \Gb \in 
X_{\Lbs}\}$. They are non-empty, by the definition of $b(\Lb)$.
By FLC, there are finitely many such
clusters up to translations. This means that there exist $\Gb_1,\dots,\Gb_K\in
X_{\Lbs}$ such that for any $\Gb\in X_{\Lbs}$ there are {\em unique}
$n=n(\Gb)\le K$ and
$u = u(\Gb)\in \R^d$ satisfying
$$
B_R(0) \cap \Gb = - u + (B_R(0) \cap \Gb_n).
$$
For $j=1,\ldots,K$ let 
$$
W_j = \{u(\Gb):\ \Gb \in X_{\Lbs}\ \mbox{such that}\ n(\Gb)=j\}.
$$
By construction, $X_{\Lbs} = \bigcup_{j=1}^K X_{\Pbs_j,W_j}$, and this is
a disjoint union. 

Next we show that the sets $W_j$ are sufficiently ``nice,''
so that they can be obtained from a finite number of closed balls using
operations of complementation, intersection, and union. 

Let $b=b(\Lb)$ and fix $j$.
Since every ball of radius $b/2$ contains a point in $\supp(\Lb)$,
we have that $W_j \subset B_{b}(0)$. Indeed, shifting a
cluster of points in $B_R(0)$ by more than $b$ would move at least one
point out of $B_R(0)$. Let $\Pb_j := B_R(0)\cap \Gb_j$.
The set $W_j$ consists of vectors
$u$ such that $-u + \Pb_j$ is a $B_R(0)$-cluster for some Delone multiset in
$X_{\Lbs}$. Thus $u\in W_j$ if and only if the following two conditions are
met. The first condition is that for each $x\in \supp(\Pb_j)$, we
have $-u + x\in B_R(0)$. The second condition is that no points of $\Gb_j$ 
outside of $B_R(0)$ move inside after the translation by $-u$.
Since $W_j \subset B_{b}(0)$, only the points in
$B_{R+b}(0)$ have a chance of moving into $B_R(0)$. Thus we need to
consider the $B_{R+b}(0)$ extensions of $\Pb_j$. 
By FLC, in the space $X_{\Lbs}$ there are finitely many
$B_{R+b}(0)$-clusters that extend the cluster $\Pb_j$.
Denote these clusters by $\Qb_1,\ldots,\Qb_L$. Summarizing this discussion
we obtain
$$
W_j = \bigcap_{x\in \suppm(\Pbs_j)}\!\!\!\!\!\!
(- B_R(0)+x) \ \ 
\cap \ \ \bigcup_{i\le L} \left[ \bigcap_{x\in \suppm(\Qbs_i) \setminus B_R(0)} 
\!\!\!\!\!\!\!\!\!\!\!\!(- (\R^d \setminus B_R(0))+x)\right].
$$
This implies that $W_j$ is a  Borel set, with $\Vol(\bd W_j)=0$.

It remains to partition each $W_j$ such that 
$W_j = \bigcup_{k=1}^{n_j} V_{jk}$, where $\diam(V_{jk}) \le \delta$, 
$0 < \delta < \eta(\Lb)$.
To this end, consider, for example, a decomposition of the cube
$[-b,b]^d$ into a
disjoint union of (half-open and closed) 
grid boxes of diameter less than $\delta < \eta(\Lb)$. Let
$\Qk$ denote the (finite) collection of all these grid boxes. 
Then 
\[ W_j = \bigcup_{D \in \Qk}(W_j \cap D) = \bigcup_{k=1}^{n_j} V_{jk} \,,
\]
where $V_{jk}$'s are disjoint and $\Vol(\partial V_{jk}) = 0$.
Note that the union $X_{\Pbs_j,W_j} = \bigcup_{k=1}^{n_j} X_{\Pbs_j,V_{jk}}$
is disjoint, from the definition of $W_j$ and $\diam(V_{jk}) < 
\eta(\Lb)$ for all $k \le n_j$.
So the lemma is proved.
\qed

\medskip

For a cluster $\Pb$ and a bounded set $A\subset \R^d$ denote
$$
L_{\Pbs}(A) = \sharp\{x\in \R^d:\ x+\Pb \subset A\cap \Lb\},
$$
where $\sharp$ means the cardinality.
In plain language, $L_{\Pbs}(A)$ is the number of translates of $\Pb$ contained
in $A$, which is clearly finite.

For a bounded set $F \subset \R^d$ and $r > 0$, let
\[
\begin{array}{l}
F^{+r} := \{x \in \R^d:\,\dist(x,F) \le r\},\\
F^{-r} := \{x \in F:\, \dist(x,\partial F) \ge r\} \supset F \setminus
(\partial F)^{+r}.
\end{array}
\]

A {\em van Hove sequence} for $\R^d$ is a sequence 
$\mathcal{F}=\{F_n\}_{n \ge 1}$ of bounded measurable subsets of 
$\R^d$ satisfying
\be \label{Hove}
\lim_{n\to\infty} \Vol((\partial F_n)^{+r})/\Vol(F_n) = 0,~
\mbox{for all}~ r>0.
\ee

\begin{defi} \label{def-ucf}
Let $\{F_n\}_{n \ge 1}$ be a van Hove sequence.
The Delone multiset $\Lb$ has {\em uniform cluster frequencies} (UCF)
(relative to $\{F_n\}_{n \ge 1}$) if for any cluster $\Pb$, the limit
$$
\freq(\Pb,\Lb) = \lim_{n\to \infty} \frac{L_{\Pbs}(x+F_n)}{\Vol(F_n)} \ge 0
$$
exists uniformly in $x\in \R^d$.
\end{defi}

Recall that a topological dynamical system is {\em uniquely ergodic} 
if there is a unique invariant probability measure. 

\begin{theorem} \label{thm-erg1}
Let $\Lb$ be a Delone multiset with FLC and $\{F_n\}_{n \ge 1}$ be a van Hove sequence.
The system  $(X_{\Lbs},\R^d)$ is uniquely ergodic if and only if
for all continuous functions $f : X_{\Lbs} \to \C ~(f\in \Cant(X_{\Lbs}))$,
\be \label{eq-erg1}
(I_n)(\Gb,f) := \frac{1}{\Vol(F_n)} \int_{F_n} f(-g+\Gb)\,dg \to\
\mbox{const},\ \ n\to \infty,
\ee
uniformly in $\Gb \in X_{\Lbs}$, with the constant depending on $f$.
\end{theorem}

This is a standard fact (both directions, see e.g. \cite[Th.\ 6.19]{Wal},
 \cite[(5.15)]{DGS}, or \cite[Th.\ IV.13]{Queff} 
for the case of $\Z$-actions); we include a (well-known)
elementary proof of the needed direction for the reader's convenience.

\medskip

\noindent {\em Proof of sufficiency in Theorem~\ref{thm-erg1}.}  
For any invariant measure $\mu$, exchanging the order of integration yields
$$
\int_{X_{\Lbs}} I_n(\Gb,f)\,d\mu(\Gb) = \int_{X_{\Lbs}} f\,d\mu,
$$
so by the Dominated Convergence Theorem, the constant in (\ref{eq-erg1}) is
$\int_{X_{\Lbs}} f\,d\mu$. If there is another invariant measure $\nu$, then
$\int_{X_{\Lbs}} f\,d\mu = \int_{X_{\Lbs}} f\,d\nu$ for all $f\in 
\Cant(X_{\Lbs})$, hence $\mu=\nu$. \qed

\medskip

Now we prove that FLC and UCF
imply unique ergodicity of the system $(X_{\Lbs},\R^d)$. 
This is also a standard fact, 
see e.g. \cite[Cor.\ IV.14(a)]{Queff} for the case of $\Z$-actions.

\begin{theorem} \label{thm-unerg} 
Let $\Lb$ be a Delone multiset with FLC. 
Then the dynamical system $(X_{\Lbs},\R^d)$ is uniquely ergodic
if and only if $\Lb$ has UCF.
\end{theorem}

\noindent {\em Proof.}
Let $X_{\Pbs,V}$ be a cylinder set with $\diam(V) \le \ \eta(\Lb)$ and $f$ 
be the characteristic function of $X_{\Pbs,V}$. 
Then we have by the definition of the cylinder set:
\be
J_n(h,f) & := & \int_{F_n} f(-x-h+\Lb)\,dx \nonumber \\
& = & \Vol\{x\in F_n:\ -x-h+\Lb\in X_{\Pbs,V}\}
\nonumber \\
& = & \Vol\{x\in h+F_n:-y+\Pb \subset -x+\Lb\ \mbox{for some}\ y\in V\}
\nonumber \\
& = & \Vol\left[\bigcup_\nu ((h+F_n) \cap (x_\nu+V))\right] \nonumber
\ee
where $x_\nu$ are all the vectors such that $x_{\nu}+\Pb \subset \Lb$.
It is clear that the distance between any two vectors $x_\nu$ is at least
$\eta(\Lb)$, so the sets $x_\nu+V$ are disjoint. Let
$$
r = \max\{|y|:\ y\in V\} + \max\{|x|:\ x\in \supp(\Pb)\}.
$$
Then
\be \label{eq-star}
\Vol(V) L_{\Pbs}(h+F_n^{-r}) \le J_n(h,f) \le \Vol(V) L_{\Pbs}(h+F_n^{+r}).
\ee
Note that
$$
L_{\Pbs}(h+F_n^{+r})-L_{\Pbs}(h+F_n^{-r}) \le L_{\Pbs}(h+{\bd F_n}^{+2r})
\le \frac{\Vol({\bd F_n}^{+2r})}{\Vol(B_{\frac{\eta({\Lbt})}{2}})}.
$$ 
So 
\be \label{uni-ucf}
\lim_{n \to \infty} \left(\frac{J_n(h,f)}{\Vol(F_n)} - \frac{\Vol(V) \cdot 
L_{\Pbs}(h + F_n)}{\Vol(F_n)} \right) = 0 \ \ 
\mbox{uniformly in}\ h \in \R^d.
\ee

If $(X_{\Lbs}, \R^d)$ is uniquely ergodic, 
$$\lim_{n \to \infty} \frac{J_n(h,f')}{\Vol(F_n)} \ \ \mbox{exists uniformly 
in} \ h \in \R^d $$ 
for continuous functions $f'$ approximating the characteristic 
function $f$ of the cylinder set. Thus for any cluster $\Pb$, 
$$\lim_{n \to \infty} \frac{L_{\Pbs}(h + F_n)}{\Vol(F_n)} \ \ 
\mbox{exists uniformly in} \ h \in \R^d,$$ i.e. $\Lb$ has UCF.

On the other hand, we assume that $\Lb$ has UCF.
By Lemma \ref{lem-tech}, $f \in \Cant(X_{\Lbs})$ can be approximated 
in the supremum norm by
linear combinations of characteristic functions of cylinder sets
$X_{\Pbs,V}$. Thus, it is enough to check
(\ref{eq-erg1}) for $f$ the characteristic function of $X_{\Pbs,V}$ with
$\diam(V)<\eta(\Lb)$.
We can see in the above (\ref{uni-ucf}) that (\ref{eq-erg1}) holds
for all $-h+\Lb$ uniformly in $h \in \R^d$ under the assumption that 
$\Lb$ has UCF. 
Then we can approximate the orbit of $\Gb \in X_{\Lbs}$ on $F_n$ by $-h_n+\Lb$ 
as close as we want, since the orbit 
$\{-h+\Lb:h \in \R^d\}$ is dense in $X_{\Lbs}$ by the definition of 
$X_{\Lbs}$.
So we compute all those integrals (\ref{eq-erg1}) of $-h_n+\Lb$ over $F_n$ 
and use the fact that independent of $h_n$ they are going to a constant.
Since each of these is uniformly close to $(I_n)(\Gb,f)$ in (\ref{eq-erg1}),
we get that $(I_n)(\Gb,f)$ too goes to a constant.
Therefore $(X_{\Lbs}, \R^d)$ is uniquely ergodic.
\qed

\medskip

Denote by $\mu$ the unique invariant probability measure on $X_{\Lbs}$.
As already mentioned,
the constant in (\ref{eq-erg1}) must be $\int_{X_{\Lbs}} f\,d\mu$. Thus, the
proof of unique ergodicity yields the following result.

\begin{cor} \label{cor-meas}
Let $\Lb$ be a Delone multiset with FLC and UCF.
Then for any $\Lb$-cluster $\Pb$ and
any Borel set $V$ with $\diam(V) < \eta(\Lb)$, we have
$$
\mu(X_{\Pbs,\tiny{V}}) = \Vol(V) \cdot\freq(\Pb,\Lb).
$$
\end{cor}

\section{Pure-pointedness and Diffraction}

\subsection{Dynamical spectrum and diffraction spectrum}

Suppose that $\Lb = (\Lam_i)_{i\le m}$ is a Delone multiset with FLC and UCF.
There are two notions of pure pointedness that appear in this context. 
Although they are defined very differently, they are in fact equivalent.  
Given a translation-bounded measure $\nu$ on $\R^d$, let $\gamma(\nu)$ denote
its autocorrelation (assuming it is unique), that is, the vague limit
\be \label{eq-auto1}
\gamma(\nu) = \lim_{n\to \infty} \frac{1}{\Vol(F_n)} \left(
\nu|_{F_n} \ast \widetilde{\nu}|_{F_n} \right),
\ee
where $\{F_n\}_{n \ge 1}$ is a van Hove sequence
\footnote{Recall that if $f$ is a function in $\R^d$, 
then $\tilde{f}$ is defined by 
$\tilde{f}(x) = \overline{f(-x)}$. If $\mu$ is a measure, $\tilde{\mu}$ is
defined by $\tilde{\mu}(f) = \overline{\mu(\tilde{f})}$ for all $f \in 
\Cant_0(\R^d)$. In particular for $\nu$ in (\ref{eq-auto2}), 
$\tilde{\nu} = \sum_{i \le m} \overline{a_i}\delta_{-\Lam_i}$.}.
In particular, for the Delone multiset $\Lb$ we see that the autocorrelation 
is unique for any measure of the form
\be \label{eq-auto2}
\nu = \sum_{i \le m} a_i \delta_{\Lambda_i},\ \ \ \mbox{where}\ \ 
\delta_{\Lambda_i} = \sum_{x\in \Lambda_i} \delta_x ~\mbox{and}~ 
a_i \in\C\,. 
\ee
Indeed, a simple computation shows
\be \label{eq-autocorr}
\gamma(\nu) = \sum_{i,j=1}^m a_i \overline{a}_j
\sum_{y \in \Lambda_i,z \in \Lambda_j} \freq((y,z), \Lb) \delta_{y-z}.
\ee
Here $(y,z)$ stands for a cluster
consisting of two points $y\in \Lambda_i, z\in \Lambda_j$.
The measure $\gamma(\nu)$ is positive definite, so by Bochner's Theorem the
Fourier transform $\widehat{\gamma(\nu)}$
is a positive measure on $\R^d$, called the {\em diffraction measure} for
$\nu$.  We say  that the measure $\nu$ 
has {\em pure point diffraction spectrum} if
$\widehat{\gamma(\nu)}$ is a pure point or discrete measure
\footnote{We also say that $\Lam_i$(resp \Lb) has pure point diffraction 
spectrum if $\widehat{\gamma(\delta_{\Lam_i})}$(resp each 
$\widehat{\gamma(\delta_{\Lam_i})}, i=1,\dots,m$) is a pure point measure.}.

\medskip

On the other hand, we also have the measure-preserving system 
$(X_{\Lbs}, \mu,\R^d)$ associated with $\Lb$. Consider the
associated group of unitary operators $\{U_x\}_{x\in \R^d}$ on
$L^2(X_{\Lbs},\mu)$:
$$
U_x f(\Lb') = f(-x + \Lb').
$$
Every $f\in L^2(X_{\Lbs},\mu)$ defines a function on $\R^d$ by
$x\mapsto (U_x f,f)$. 
This function is positive definite on $\R^d$, so its
Fourier transform is a positive measure $\sigma_f$ on $\R^d$ called the
{\em spectral measure} corresponding to $f$.
We say that the Delone multiset $\Lb$ has {\em pure point dynamical spectrum}
if $\sigma_f$ is pure point for every $f\in L^2(X_{\Lbs},\mu)$.
We recall that $f \in L^2(X_{\Lbs},\mu)$
is an eigenfunction for the $\R^d$-action if for some
$\alpha =(\alpha_1,\ldots,\alpha_d) \in \R^d$,
$$
U_x f = e^{2 \pi i x\cdot \alpha} f,
\ \ \ \mbox{for all}\ \ x\in \R^d,
$$
where $\cdot$ is the standard inner product on $\R^d$.

\begin{theorem} \label{thm-pp}
$\sigma_f$ is pure point for every $f\in L^2(X_{\Lbs},\mu)$
if and only if the eigenfunctions for the $\R^d$-action span a dense
subspace of $L^2(X_{\Lbs},\mu)$.
\end{theorem} 

This is a straightforward consequence of the Spectral Theorem, see e.g.\
Theorem 7.27 and \S 7.6 in \cite{Weidmann} for the case $d=1$.
The Spectral Theorem for unitary representations of
arbitrary locally compact Abelian groups, including $\R^d$, is discussed
in \cite[\S 6]{Mackey}.

\subsection{An equivalence theorem for pure pointedness}

In this section we prove the following theorem.

\begin{theorem} \label{thm-equiv} Suppose that a Delone multiset 
$\Lb$ has FLC and UCF. Then the following are equivalent:
\begin{itemize} 
\item[{\rm (i)}] $\Lb$ has pure point dynamical spectrum;
\item[{\rm (ii)}] 
The measure $\nu = \sum_{i \le m} a_i \delta_{\Lambda_i}$ has pure point
diffraction spectrum, for any choice of complex numbers $(a_i)_{i\le m}$;
\item[{\rm (iii)}] 
The measures $\delta_{\Lambda_i}$ have pure point diffraction spectrum,
for $i\le m$.
\end{itemize}
\end{theorem}

\medskip

\noindent
The theorem is proved after a sequence of auxiliary lemmas. Fix complex 
numbers $(a_i)_{i \le m}$ and let $\nu = \sum_{i\le m} a_i \delta_{\Lambda_i}$. 
For $\Lb' = (\Lam'_i)_{i\le m} \in X_{\Lbs}$ let
$$
\nubp = \sum_{i\le m} a_i \delta_{\Lam'_i},
$$
so that $\nu = \nub$. To relate the autocorrelation of $\nu$ to
spectral measures we need to do some ``smoothing.''
Let $\om \in \Cant_0(\R^d)$
(that is, $\om$ is continuous and has compact support). Denote
$$
\rho_{\om,{\Lbs}'} := \om \ast \nubp
$$
and let
$$
f_\om(\Lb') := \rho_{\om,{\Lbs}'}(0)\ \ \ \mbox{for}\ \ \Lb' \in X_{\Lbs}.
$$

\begin{lemma} \label{lem-cont}
$f_\om \in \Cant(X_{\Lbs}).$
\end{lemma}

\noindent {\em Proof.} We have
$$
f_\om(\Lb') = \int \om(-x)\,d\nubp(x) = \sum_{i\le m} a_i 
\sum_{x\in -\suppm(\om) \cap \Lambda_i'} \om(-x).
$$
The continuity of $f_\om$ follows from the continuity of $\om$ and the
definition of topology on $X_{\Lbs}$. \qed

\medskip

Denote by $\gam_{\om,{\Lbs}}$ the autocorrelation of $\rho_{\om,{\Lbs}}$.
Since under our assumptions there is a unique autocorrelation measure
$\gam=\gamma(\nu)$, see (\ref{eq-auto1}) and (\ref{eq-auto2}), we have
$$
\gam_{\om,{\Lbs}} = (\om \ast \widetilde{\om}) \ast \gamma.
$$

\begin{lemma} \label{lem-spectr} {\rm (\cite{Dwor}, see also  \cite{Hof})} 
$$\sigma_{f_\om} = \widehat{\gamma_{\om,{\Lbs}}}.$$
\end{lemma}

\noindent {\em Proof.} 
We provide a proof for completeness, following \cite{Hof}.
By definition,
$$
f_\om(-x + \Lb) = \rho_{\om,{\Lbs}}(x).
$$
Therefore,
\be \label{Dwor-argu}
\gamma_{\om,{\Lbs}}(x) & = & \lim_{n\to\infty} \frac{1}{\Vol(F_n)}
\int_{F_n}\rho_{\om,{\Lbs}}(x+y)\overline{\rho_{\om,{\Lbs}}(y)}
\,dy \nonumber \\
& = & \lim_{n\to\infty} \frac{1}{\Vol(F_n)} \int_{F_n}
f_{\om}(-x-y+\Lb)\overline{f_\om(-y+\Lb)}\,dy \nonumber \\
& = & \int_{X_{\Lbs}} f_\om(-x+\Lb')\overline{f_\om(\Lb')}\,d\mu(\Lb')
\nonumber \\
& = & (U_x f_\om, f_\om)\,,
\ee
where $\{F_n\}_{n \ge 1}$ is a van Hove sequence.
Here the third equality is the main step; it follows from unique ergodicity
and the continuity of $f_\om$, see Theorem~\ref{thm-erg1}.
Thus,
$$
\widehat{\gamma_{\om,{\Lbs}}} = \widehat{(U_{(\cdot)} f_\om, f_\om)} = \sigma_{f_\om},
$$
and the proof is finished. \qed

The introduction of the function $f_{\om}$ and the series of equations 
(\ref{Dwor-argu}) is often called Dworkin's argument.

\medskip

Fix $\e$ with $0 < \e < \frac{1}{b(\Lb)}$. Consider all the clusters of 
diameter $\le 1/\e$ in $\Gb \in
X_{\Lbs}$. There are finitely many such clusters up to translation, by FLC.
Thus, there exists $0 < \theta_1=\theta_1(\e) <1$ 
such that if $\Pb,\Pb'$ are two such clusters, then
\be \label{eq-trans}
\rho_H(\Pb,\Pb')\le \theta_1\ \Rightarrow\ \Pb=-x+\Pb'\ \ \ \mbox{for some}\ 
x\in \R^d.\
\ee
Here 
$$\rho_H(\Pb,\Pb') = {\mmax}\{\rho_H(P_i,P'_i):\  i\le m\},
$$
 where 
$$
\rho_H(P_i,P'_i) = \left\{ \begin{array}{l} 
      {\mmax}\{ \dist(x,P'_i), \dist(y,P_i) : x \in P_i, y \in P'_i\},
      \ \ \ \mbox{if}~ P_i,P'_i \neq \emptyset;\\
        1, \ \ \ \mbox{if}~  P_i = \emptyset~ \mbox{and}~ P'_i \neq \emptyset 
                         \ (\mbox{or vice versa}),
                             \end{array} 
                   \right.
$$
with $\Pb=(P_i)_{i\le m}$ and $\Pb' = (P'_i)_{i\le m}$\,.

\smallskip

Let 
\be \label{def-theta}
\theta=\theta(\e) := {\mmin}\{\e,\theta_1,\eta(\Lb)\}
\ee
and 
$$
f_{i,\om}(\Lb') = (\om\ast\delta_{\scriptsize {\Lam_i'}})(0)\ \ \ \mbox{for}\ \ 
\Lb' = (\Lambda_i')_{i\le m} \in X_{\Lbs}.
$$
Denote by ${\bf E}_i$ the cluster consisting of a single point of type $i$ at the
origin; formally, 
$$
{\bf E}_i = (\es,\ldots,\es,\,\underbrace{\{0\}}_{i},\,\es,\ldots,\es).
$$
Let $\cheiv$ be the characteristic function for the
cylinder set $X_{{\scriptsize {\bf E}_i},V}$.

\begin{lemma} \label{lem-vspom}
Let $V\subset \R^d$ be a bounded set with 
$\diam(V) < \theta$, where $\theta$ is defined by (\ref{def-theta}), and 
$0 < \zeta < \theta/2$. 
Let $\om\in \Cant_0(\R^d)$ be such that
$$
\left\{ \begin{array}{ll} \om(x) = 1, & x\in V^{-\zeta};\\
                          \om(x) = 0, & x\in \R^d \setminus V;\\
                          0\le\om(x)\le 1, & x\in V\setminus V^{-\zeta}.
        \end{array} \right.
$$
Then
$$
\|f_{i,\om} - \cheiv\|_2^2 \le \freq({\bf E}_i,\Lb) \cdot\Vol((\bd V)^{+\zeta}).
$$
\end{lemma}

\noindent {\em Proof.} We have by the definition of ${\bf E}_i$ and Definition
\ref{def-cyl}:
$$
\cheiv(\Lb') = \left\{ \begin{array}{ll} 1, & \mbox{if}\ \Lambda_i'\cap
(-V) \ne \es;\\
0, & \mbox{otherwise}, \end{array} \right.
\ \ \ \mbox{where}\ \Lb' \in X_{\Lbs}.
$$
On the other hand, since $\om$ is supported in $V$ and there is at most one
point of $\Lambda_i'$ in $V$,
$$
f_{i,\om}(\Lb') = \int \om(-x) \,d\delta_{{\Lam}'_i}(x) = 
\left\{ \begin{array}{ll} \om(-x), & \mbox{if\ } \exists\,x\in \Lambda_i' \cap
(-V);\\  0, & \mbox{otherwise}.\end{array} \right.
$$
It follows that 
$$
f_{i,\om}(\Lb') - \cheiv(\Lb') = 0\ \ \ \mbox{if}\ \ \ 
\Lambda_i' \cap (-V^{-\zeta}) \ne \es.
$$
Thus,
\be
\|f_{i,\om} - \cheiv\|_2^2 & \le &
\int_{X_{{\bf E}_i, V\setminus V^{-\zeta}}} |f_{i,\om}(\Lb') - 1|^2\,d\mu(\Lb') 
\nonumber \\
& \le &  \mu(X_{{\bf E}_i, V\setminus V^{-\zeta}}) \nonumber \\
& =   & \freq({\bf E}_i,\Lb)\cdot \Vol(V\setminus V^{-\zeta}) \nonumber \\
& \le & \freq({\bf E}_i,\Lb)\cdot \Vol((\bd V)^{+\zeta}),\nonumber
\ee
as desired. \qed

\medskip

\begin{lemma} \label{lem-prod}
Let $\Pb = (P_i)_{i\le m} = B_{1/{\e}}(0) \cap \Gb$ with 
$\Gb \in X_{\Lbs}$, and $\diam(V) < \theta$, where $\theta$ is defined by 
(\ref{def-theta}). 
We have
$$
\chipv = \prod_{i\le m} \prod_{x\in P_i} \chi_{_{x+{\scriptsize {\bf E}_i},V~.}}
$$
\end{lemma}

\noindent {\em Proof.} We just have to prove that
$$
X_{\Pbs,V} = \bigcap_{i\le m} \bigcap_{x \in P_i} 
X_{x+{\scriptsize {\bf E}_i},V}.
$$
A Delone multiset $\Gb$ is in the left-hand side whenever $-v+\Pb\subset 
\Gb$ for some $v\in V$. 
A Delone multiset $\Gb$ is in the right-hand side whenever 
for each $i\le m$ and each $x\in P_i$ there is a vector $v({\bf x})\in V$ 
such that $-v({\bf x})+{\bf x} \subset \Gb$, where ${\bf x}=(\emptyset,\dots,
\emptyset,\underbrace{\{x\}}_{i},\emptyset,\dots,\emptyset)$ stands for a single element 
cluster. Thus, ``$\subset$'' is trivial.
The inclusion ``$\supset$''
follows from the fact that $\diam(V) < \theta$,
see (\ref{def-theta}) and (\ref{eq-trans}). \qed

\medskip

Denote by $\Hk_{pp}$ the closed linear span in $L^2(X_{\Lbs},\mu)$
of the eigenfunctions for the dynamical system $(X_{\Lbs},\mu,\R^d)$.
The following lemma is certainly standard, but since we do not know a 
ready reference, a short proof is provided.

\begin{lemma} \label{lem-prod2}
If $\phi$ and $\psi$ are both in $L^\infty(X_{\Lbs},\mu) \cap \Hk_{pp}$, then
their product $\phi\psi$ is in $L^\infty(X_{\Lbs},\mu) \cap \Hk_{pp}$ as well.
\end{lemma}

\noindent {\em Proof.} Fix arbitrary $\epsilon>0$. Since $\phi \in \Hk_{pp}$, we can 
find a finite linear combination 
of eigenfunctions $\widetilde{\phi} = \sum a_i f_i$ such that
$$
\| \phi - \widetilde{\phi} \|_2 < \frac{\epsilon}{\|\psi\|_\infty}\,.
$$
Since the dynamical system is ergodic, the eigenfunctions have constant
modulus, hence $\widetilde{\phi}\in L^\infty$.
Thus, we can find another finite linear combination 
of eigenfunctions $\widetilde{\psi} = \sum b_j f_j$ such that
$$
\| \psi - \widetilde{\psi} \|_2 < \frac{\epsilon}{\|\widetilde{\phi}\|_\infty}\,.
$$
Then
\be 
\| \phi \psi - \widetilde{\phi} \widetilde{\psi} \|_2 & \le & 
\|\widetilde{\phi}(\psi-\widetilde{\psi})\|_2 + \
\|(\phi - \widetilde{\phi})\psi\|_2 \nonumber \\
& \le &  \|\widetilde{\phi}\|_\infty \| \psi - \widetilde{\psi} \|_2
+ \|\psi\|_\infty \| \phi - \widetilde{\phi} \|_2 \nonumber \\
& \le & 2\epsilon. \nonumber
\ee
It remains to note that $\widetilde{\phi} \widetilde{\psi}\in \Hk_{pp}$
since the product of eigenfunctions for a dynamical system  is an eigenfunction.
Since $\epsilon$ is arbitrarily small, $\phi\psi\in \Hk_{pp}$,
and the lemma is proved.  \qed

\medskip

\noindent {\em Proof of Theorem~\ref{thm-equiv}.} 
(i) $\Rightarrow$ (ii)
This is essentially proved by Dworkin in \cite{Dwor}, see also  
\cite{Hof} and \cite{BM}. By Lemma~\ref{lem-spectr}, pure point dynamical
spectrum implies that $\widehat{\gamma_{\om,{\Lbs}}}$ is pure point for any
$\om \in \Cant_0(\R^d)$. Note that
\be \label{eq-auto3}
\widehat{\gamma_{\om,{\Lbs}}} = |\widehat{\om}|^2 \widehat{\gamma}.
\ee
Choosing a sequence $\om_n \in \Cant_0(\R^d)$ converging to the delta measure 
 $\delta_0$ in the
vague topology, we can conclude that $\widehat{\gamma}$ is pure point
as well, as desired. (This approximation step requires some care; it is
explained in detail in \cite{BM}.)

\smallskip

(ii) $\Rightarrow$ (iii) obvious.

\smallskip

(iii) $\Rightarrow$ (i) 
This is relatively new, although it
is largely a generalization of Queffelec \cite[Prop.\ IV.21]{Queff}.

We are given that $\delta_{\Lambda_i}$ has pure point diffraction spectrum, that is,
$\widehat{\gam_i}:=\widehat{\gamma(\delta_{\Lambda_i})}$ is pure point, for all $i\le m$.
In view of (\ref{eq-auto3}) and Lemma~\ref{lem-spectr}, we obtain that
$\sigma_{f_{i,\om}}$ is pure point for all $i\le m$ and all $\om \in
\Cant_0(\R^d)$. So $f_{i,\om} \in \Hk_{pp}$ for all $i\le m$ and all $\om \in
\Cant_0(\R^d)$.
Fix $\e > 0$ and let $V$ be a bounded set with $\diam(V) < \theta=\theta(\e)$,
where $\theta$ is defined by (\ref{def-theta}), and $\Vol(\bd V) = 0$. 
Find $\om\in \Cant_0(\R^d)$ as in Lemma~\ref{lem-vspom}. 
Since $\Vol((\bd V)^{+\zeta}) \to 
\Vol(\bd V) = 0$ in Lemma~\ref{lem-vspom}, as $\zeta \to 0$,
we obtain that $\cheiv \in \Hk_{pp}$.
Therefore, also $U_x\cheiv = \chi_{_{x+{\scriptsize {\bf E}_i},V}}\in 
\Hk_{pp}$.
Then it follows from Lemma \ref{lem-prod} and 
Lemma \ref{lem-prod2} that $\chipv \in \Hk_{pp}$ where $\Pb = B_{1/\e}(0)\cap \Gb$ for any $\Gb \in X_{\Lbs}$, $\diam(V) < \theta$, and $\Vol(\bd V) = 0$. 

Our goal is to show that $\Hk_{pp} = L^2(X_{\Lbs},\mu)$.
Since $(X_{\Lbs},\mu)$ is a regular measure space, $\Cant(X_{\Lbs})$ is dense
in $L^2(X_{\Lbs},\mu)$. Thus, it is enough to show that all continuous functions
on $X_{\Lbs}$ belong to $\Hk_{pp}$. Fix $f\in \Cant(X_{\Lbs})$.
Using the decomposition $X_{\Lbs} = \bigcup_{j=1}^N X_{\Pbs_j,V_j}$ 
from Lemma~\ref{lem-tech} we can approximate
$f$ by linear combinations of characteristic functions of cylinder sets 
$X_{\Pbs_j,V_j}$.
So it suffices to show that these characteristic functions are in $\Hk_{pp}$,
which was proved above.
This concludes the proof of Theorem \ref{thm-equiv}. \qed

\section{Concluding remarks: what if the UCF fails?}

Here we present a version of the main theorem for Delone sets which do not 
necessarily have uniform cluster frequencies.
For this we must assume that in addition to the van Hove property 
(\ref{Hove}) our averaging sequence $\{F_n\}$ is a sequence of compact 
neighbourhoods of $0$ with the following properties:
\be  \label{add-van Hove}
\begin{array}{l} 
\mbox{(i) $\cup F_n = \R^d$} \\
\mbox{(ii) $\exists$  $K \ge 1$ so that 
$\Vol(F_n - F_n) \le K \cdot \Vol(F_n)$ for all $n$}.
\end{array}
\ee
Let $\Lb\subset \R^d$ be a Delone set with FLC. Consider the topological
dynamical system $(X_{\Lbs},\R^d)$ and an ergodic invariant Borel
probability measure $\mu$ (such measures always exist). 
The ergodic measure $\mu$ will be fixed throughout the section.

\begin{theorem} \label{th-noucf}
Suppose that a Delone multiset $\Lb$ has FLC. 
Then the following are equivalent:
\begin{itemize} 
\item[{\rm (i)}] The measure-preserving dynamical system $(X_{\Lbs},\mu,\R^d)$
has pure point spectrum;
\item[{\rm (ii)}] For $\mu$-a.e.\ $\Gb\in X_{\Lbs}$, 
the measure $\nu = \sum_{i \le m} a_i \delta_{\Gamma_i}$ has pure point
diffraction spectrum, for any choice of complex numbers $(a_i)_{i\le m}$;
\item[{\rm (iii)}] For $\mu$-a.e.\ $\Gb\in X_{\Lbs}$,
the measures $\delta_{\Gamma_i}$ have pure point diffraction spectrum,
for $i\le m$.
\end{itemize}
\end{theorem}

In fact, this formulation is closer to the work of Dworkin \cite{Dwor}
who did not assume unique ergodicity. The proof is similar to that of 
Theorem~\ref{thm-equiv}, except that we have to use the Pointwise
Ergodic Theorem instead of the uniform convergence of averages
in the uniquely ergodic case (\ref{eq-erg1}).

\begin{theorem} \label{thm-pet}
(Pointwise Ergodic Theorem for $\R^d$-actions (see, e.g. \cite{Cha}))
For any $f\in L^1(X_{\Lbs},\mu)$,
\be \label{th-pet}
\frac{1}{\Vol(F_n)} \int_{F_n} f(-x+\Gb)\,dx \to \int f(\Lb')\,d\mu(\Lb'),\ \ 
\mbox{as}\ n\to \infty,
\ee
for $\mu$-a.e.\ $\Gb \in X_{\Lbs}$.
\end{theorem}

For a cluster $\Pb\subset \Lb$, a bounded set $A\subset \R^d$, and a 
Delone set $\Gb\in X_{\Lbs}$, denote 
$$L_{\Pbs} (A,\Gb) = \# \{x\in \R^d:\ x + \Pb \subset A \cap \Gb\}.$$

\begin{lemma} \label{lem-nounif}
For $\mu$-a.e.\ $\Gb \in X_{\Lbs}$ and for any cluster $\Pb \subset
\Lb$, 
\be \label{N1}
\freq'(\Pb,\Gb):= \lim_{n \to \infty} \frac{L_{\Pbs}(F_n,\Gb)}{\Vol(F_n)},
\ee 
exists for $\mu$-a.e.\ $\Gb \in X_{\Lbs}$.
Moreover, if $\diam(V) < \eta(\Lb)$, then the cylinder set $X_{\Pbs,V}$
satisfies, for $\mu$-a.e.\ $\Gb \in X_{\Lbs}$:
\be \label{N2}
\mu(X_{\Pbs,V}) = \Vol(V) \cdot \freq'(\Pb,\Gb).
\ee
\end{lemma}

Note that we no longer can claim uniformity of the convergence with 
respect to translation of $\Gb$.

\medskip

\noindent {\em Sketch of the proof.} Fix a cluster $\Pb\subset \Lb$ and  
let $X_{\Pbs,V}$ be a cylinder set, with $\diam(V) < \eta(\Lb)$.
Applying (\ref{th-pet}) to the characteristic function of $X_{\Pbs,V}$
and arguing as in the proof of Theorem~\ref{thm-unerg}
(with $-h + \Lb$ replaced by $\Gb$), 
we obtain (\ref{N1}) and
(\ref{N2}) for $\mu$-a.e.\ $\Gb$. Since there are countably many
clusters $\Pb \subset \Lb$, we can find a set of full $\mu$-measure
on which (\ref{N1}) and (\ref{N2}) hold for all $\Pb$. \qed

\medskip

For a Delone set $\Gb=(\Gamma_i)_{i\le m}$, let 
$\nu =\sum_{i=1}^m a_i \delta_{\Gamma_i}$.
Then, for $\mu$-a.e.\ $\Gb$, the autocorrelation
$\gamma(\nu)$ exists as the vague limit of measures
$\frac{1}{\Vol(F_n)} (\nu|_{F_n}\ast\widetilde{\nu}|_{F_n})$, and 
$$
\gamma(\nu) = \sum_{i,j=1}^m a_i \overline{a}_j
\sum_{y \in \Gamma_i,z \in \Gamma_j} \freq'((y,z), \Gb) \delta_{y-z},
$$
for $\mu$-a.e.\ $\Gb$, which is the analogue of (\ref{eq-autocorr}).
Again, $\widehat{\gamma(\nu)}$ is a positive
measure, called the diffraction measure, giving the meaning to the words
``pure point diffraction spectrum'' in Theorem~\ref{th-noucf}.

\medskip

\noindent {\em Sketch of the proof of Theorem~\ref{th-noucf}.}
For $\mu$-a.e.\ $\Gb\in X_{\Lbs}$, the Pointwise Ergodic Theorem
\ref{thm-pet} holds for all functions $f \in \Ck(X_{\Lbs})$ (since the
space of continuous functions on $X_{\Lbs}$ is separable).

The $\nu_{\Lbt'}, \rho_{\om,\Lbt'}$, and $f_\om$ are defined the same
way as in Section 3. Lemma~\ref{lem-cont} applies to our
situation. Next we can show that 
\be \label{eq-new}
\sigma_{f_\om} = \widehat{\gamma_{\om,\Gbt}}
\ee
for $\mu$-a.e.\ $\Gb$.
This is proved by the same chain of equalities as in (\ref{Dwor-argu}),
except that we average over $F_n$ defined in (\ref{add-van Hove}) 
and use Theorem~\ref{thm-pet} instead of Theorem~\ref{thm-unerg}. 
Lemma~\ref{lem-vspom} goes through, after we replace $\freq({\bf E}_i,\Lb)$ 
by $\freq'({\bf E}_i,\Gb)$, for $\mu$-a.e.\ $\Gb$.
There are no changes in Lemmas~\ref{lem-prod} and \ref{lem-prod2},
since we did not use UCF or unique ergodicity in them. The proof of 
Theorem~\ref{th-noucf} now follows the scheme of the proof of
Theorem~\ref{thm-equiv}. We only need to replace $\Lb$ by $\mu$-a.e.\ 
$\Gb$, for which hold all the ``typical'' properties discussed above.
\qed

\end{document}